\documentclass[12pt]{amsart}
\usepackage{mathrsfs}
\usepackage{amsmath,amssymb,amsfonts,latexsym,txfonts}
\usepackage{tikz}
\newlength{\originalbase}
\newcommand{\spacing}[1]{\setlength{\baselineskip}{#1\originalbase}}

\newtheorem{theorem}{Theorem}[section]
\newtheorem{proposition}[theorem]{Proposition}

\newtheorem{corollary}[theorem]{Corollary}

\newcommand{\RR}{\mathbb{R}}
\newcommand{\uu}{{\bf u}}
\newcommand{\vv}{{\bf v}}
\newcommand{\ww}{{\bf w}}
\newcommand{\oo}{{\bf o}}

\begin{document}

\title{Tetrahedra with congruent face pairs}
\author[D.~Klain]{Daniel A. Klain} 
\address{D.~Klain, Department of Mathematical Sciences,
University of Massachusetts Lowell,
Lowell, MA 01854, USA}
\email{Daniel\_{}Klain@uml.edu}

\begin{abstract} 
If the four triangular facets of a tetrahedron 
can be partitioned into pairs having the same area, then 
the triangles in each pair must be congruent to one another.
A Heron-style formula is then derived for the volume of a tetrahedron
having this kind of symmetry.\\[1mm]
\end{abstract}

\maketitle

\vspace{4mm}

\spacing{1}

From elementary geometry we learn that two triangles are congruent if
their edges have the same three lengths. In particular, there is only
one congruence class of {\em equilateral} triangles having a given edge
length.  Said differently, any pair of equilateral triangles in the
Euclidean plane are {\em similar}, differing at most by an isometry 
and a dilation.  Meanwhile, triangles that are symmetric under a single reflection
have two congruent sides and are said to be {\em isosceles}.

The situation is more complicated in higher dimensions.  Indeed, an 
analogous characterization of $3$-dimensional tetrahedra 
already leads to $25$ different symmetry classes \cite{WD2}.
These tetrahedral symmetry classes are of special interest in organic chemistry \cite{FR1,FR2,WD1},
and conditions for tetrahedral symmetry based on the measures of dihedral angles have also
been explored \cite{WD3}.

A tetrahedron in $\RR^3$ is {\em equilateral} or {\em regular} if
all of its edges have the same length. 
More generally, a tetrahedron is said to be 
{\em isosceles} if 
all four triangular facets are congruent
to one another, or, equivalently, if opposing (non-incident) edges 
have the same length. Isosceles tetrahedra are also known as 
{\em disphenoids} \cite[p.~15]{Coxeter}.  
It has been shown that
if all four facets of a tetrahedron $T$ have the same {\em area}, 
then $T$ must be isosceles \cite[p.~94]{GemsII}\cite{Horv,equiareal}.

Consider the following more general symmetry class of tetrahedra: 
A tetrahedron $T$ will be called {\em  reversible} if its facets
are congruent in pairs; that is, if the facets of $T$ can be labelled
$f_1, f_2, f_3, f_4$,
where $f_1 \cong f_2$ and $f_3 \cong f_4$.  

\begin{figure}[ht]
\begin{tikzpicture}[scale=1.2]
\draw[thick] (0,0)--(2,2)--(7,1)--(5,-1)--(0,0);
\draw[thick,dotted] (0,0)--(7,1);
\draw[thick] (2,2)--(5,-1);
\node [below] at (2.5,-0.5) {\small $a$};
\node [above] at (4.5,1.5) {\small $a$};
\node [left] at (0.9,1) {\small $b$};
\node [right] at (6.1,0) {\small $b$};
\node [above] at (2,0.25) {\small $d$};
\node [above] at (4.5,-0.4) {\small $c$};
\end{tikzpicture}
\caption{An reversible tetrahedron with edge lengths $a,a,b,b,c,d$.} 
\label{pn-fig}
\end{figure}
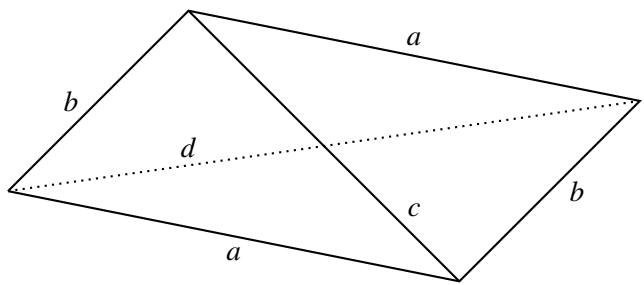

In this note we show
that, as in the isosceles case, reversible tetrahedra are characterized by the areas of their facets:  
if the four triangular facets of $T$ can be partitioned into pairs with the same area,
then those pairs consist of congruent facets.

In the final section we give an intuitive method for deriving a 
Heron-style factorization of the volume of a reversible tetrahedron 
in terms of its edge lengths.

\section{\bf Facets normals and areas determine tetrahedra}

The following proposition will allow us to exploit symmetries more easily.
\begin{proposition}
Suppose that $u_0, u_1, u_2, u_3$ are unit vectors that span $\RR^3$,
and suppose that $\alpha_0, \alpha_1, \alpha_2, \alpha_3 > 0$.  
Then there exists a tetrahedron $T$, having outward facet unit normals
$u_i$, and corresponding
facet areas $\alpha_i$, if and only if
\begin{align}
\alpha_0 u_0 + \alpha_1 u_1 + \alpha_2 u_2 + \alpha_3 u_3 = 0.
\label{mink2}
\end{align}
Moreover, this tetrahedron is unique up to translation.
\label{special}
\end{proposition}

This proposition is a very special case of the Minkowski Existence Theorem, which plays a central role
in the Brunn-Minkowski theory of convex bodies, 
and is somewhat difficult to prove \cite{Bonn,red}.
However, this special case for tetrahedra is a simple consequence of linear algebra.

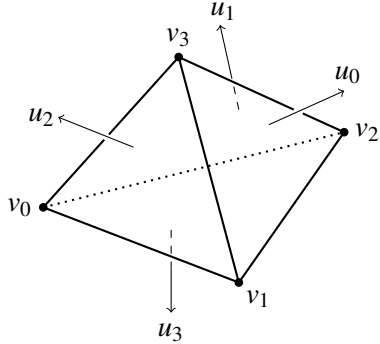
\begin{figure}[ht]
\begin{tikzpicture}[scale=1]
\draw[thick] (0,0)--(0.81,.9);
\draw[thick] (0.873,0.97)--(1.8,2);
\draw[thick] (0,0)--(2.6,-1);
\draw[thick,dotted] (0,0)--(4,1);
\draw[thick] (1.8,2)--(2.6,-1)--(4,1);
\draw[thick] (4,1)--(3.45,1.25);
\draw[thick] (3.33,1.3)--(1.8,2);
\draw[fill] (0,0) circle [radius=0.05];
\draw[fill] (1.8,2) circle [radius=0.05];
\draw[fill] (2.6,-1) circle [radius=0.05];
\draw[fill] (4,1) circle [radius=0.05];
\draw[dashed] (1.7,-0.3)--(1.7,-0.8);  
\draw[->] (1.7,-0.8)--(1.7,-1.4); \node [below] at (1.7,-1.4) {\small $u_3$};
\draw[dashed] (2.6,1.3)--(2.5,1.75);
\draw[->] (2.5,1.75)--(2.35,2.425); \node [above] at (2.4,2.35) {\small $u_1$};
\draw[->] (1.2,0.8)--(0.2,1.2); \node [left] at (0.3,1.2) {\small $u_2$};
\draw[->] (3,1.1)--(4,1.55); \node [above] at (4.05,1.5) {\small $u_0$};
\node [left] at (0,0) {\small $v_0$};
\node [above] at (1.8,2) {\small $v_3$};
\node [below right] at (2.55,-0.95) {\small $v_1$};
\node [right] at (4,1) {\small $v_2$};
\end{tikzpicture}
\caption{A tetrahedron with outward unit normals $u_i$.} 
\label{tetra-fig}
\end{figure}

\begin{proof} Let $T$ be a tetrahedron with vertices at $v_0, v_1, v_2, v_3 \in \RR^3$, where $v_0 = o$, the origin.
Let us assume the vertices are labelled so that $v_1, v_2, v_3$ have a positive (``right-handed'') orientation.

Denote by $u_0, u_1, u_2, u_3$ the outward unit normal vectors of the facets of $T$, 
where $u_i$ is associated with the facet opposite to the vertex $v_i$, 
as in Figure~\ref{tetra-fig}.    
Let $\alpha_i$ denote the area of that same $i$th facet.  Since $v_0 = o$,
we have
\begin{align} \notag
v_2 \times v_3 &= -2\alpha_1 u_1\\ \label{crosses}
v_3 \times v_1 &= -2\alpha_2 u_2\\ \notag
v_1 \times v_2 &= -2\alpha_3 u_3\\ \notag
(v_3-v_1) \times (v_2-v_1) &= -2 \alpha_0 u_0.
\end{align}
After summing both sides of these equations the identity~(\ref{mink2}) now follows.

Next, let $A$ denote a $3 \times 3$ matrix whose columns are given by the vectors $v_1, v_2, v_3$.  
Recall that these vectors are ordered so that $A$ has positive determinant.  
In this instance $\det(A) = 6V(T)$, where $V(T)$ denotes the volume of $T$. 

Let $c(A)$ denote the cofactor matrix of $A$.  Cramer's Rule asserts that
\begin{equation}
c(A)^t A = \det(A) I,
\label{cramer}
\end{equation}
where $I$ is the $3 \times 3$ identity matrix.  (See \cite[p.~30]{Artin}, for example,  
or any traditional linear algebra text.)  

Let $z_i$ denote the $i$th column of the matrix $c(A)$.  The
identity~(\ref{cramer}) asserts that $z_i \perp v_j$ for $j \neq i$. 
It follows from~(\ref{crosses}) that $z_i$ is parallel to the facet normal $u_i$, 
and that $z_i = -|z_i| u_i$, 
since $z_i \cdot v_i = \det(A) > 0$, 
while $u_i$ points out of the tetrahedron (away from the vertex $v_i$).
Meanwhile,~(\ref{cramer}) also asserts that
$$z_i \cdot v_i \; = \; \det(A)  \; = \; 6 V(T).$$
It follows that
$$-|z_i| (u_i \cdot v_i) \; = \; z_i \cdot v_i \; = \; 6V(T) 
\; = \; 2\alpha_i (-u_i \cdot v_i),$$
where the final identity follows from the base-height formula for the
volume of a cone, using the $i$th facet of $T$ as the base. Hence,
$|z_i| = 2\alpha_i $ and $z_i = -2\alpha_i u_i.$ In
other words, the facet normals $u_1, u_2, u_3$ and corresponding
facet areas $\alpha_i$ are determined by the columns $z_i$ of the cofactor
matrix $c(A)$. The remaining facet normal $u_0$ and area $\alpha_0$ is
then determined by the identity~(\ref{mink2}).  
This encoding of facet data into the cofactor
matrix allows a simple proof of both existence and uniqueness for the
tetrahedron $T$ given the data $\{u_i\}$ and $\{\alpha_i\}$.

To prove the {\em uniqueness} of $T$, note that $c(A) = \det(A) A^{-t}$, 
by Cramer's Rule~(\ref{cramer}).  
It follows that $\det(c(A)) = \det(A)^2$ and that
$$A \; = \; \det(A) c(A)^{-t} \; = \; \det(c(A))^{\frac{1}{2}} c(A)^{-t}.$$
In other words, if two matrices $A$ and $B$ with positive determinant 
have the same
cofactor matrix $c(A) = c(B)$, then $A = B$.  It follows that if two
tetrahedra $T_1$ and $T_2$ each have the origin as a vertex and share the
same facet normals and corresponding facet areas, 
then $T_1$ and $T_2$ must have the same vertices, so
that $T_1 = T_2$.

More generally, if two tetrahedra have the same facet normals and
corresponding facet areas, then they must be translates of one another.

To prove the {\em existence} of a tetrahedron 
having the given facet data,
Let $C$ denote the matrix having columns $-2\alpha_i u_i$ for
$i > 0$, ordered so that $C$ has positive determinant.  The matrix 
$$A = \det(C)^{\frac{1}{2}} C^{-t}$$
has cofactor matrix $C$. The
columns of $A$, along with the origin, yield the vertices of a tetrahedron
having facet normals $u_i$ and corresponding facet areas $\alpha_i$.
\end{proof}

\section{\bf Equal areas imply congruent faces}

We now prove that the areas of the facets alone will determine if a 
tetrahedron is reversible.
\begin{theorem}
Suppose that $T$ is a tetrahedron in $\RR^3$, and denote by
$f_1, f_2, f_3, f_4$ the triangular facets of $T$.  
If the facets of $T$ satisfy the conditions
$$Area(f_1) = Area(f_2) \;\;\; \hbox{ and } \;\;\; Area(f_3) = Area(f_4)$$
then $f_1 \cong f_2$ and $f_3 \cong f_4$.
\label{pair}
\end{theorem}

The proof of Theorem~\ref{pair} uses the method given by McMullen in \cite{equiareal}
to verify the special case in which all four facets have the same area 
(as in Corollary~\ref{regular} below).
\begin{proof}
Denote by $u_i$ the outward unit normal vector to the 
facet $f_i$ of $T$.  Suppose that
$Area(f_1) = Area(f_2) = \alpha$ and
$Area(f_3) = Area(f_4) = \beta$,
where $\alpha, \beta > 0$.  The identity~(\ref{mink2}) asserts that
$$\alpha u_1 + \alpha u_2 + \beta u_3 + \beta u_4 = 0.$$
Denote
$$w = \alpha u_1 + \alpha u_2 = -\beta u_3 - \beta u_4.$$

Let $\psi$ denote the rotation of $\RR^3$ by the angle $\pi$ around the
the axis through $w$.  Since the vectors $\alpha u_1$ and $\alpha u_2$
have the same length, the points $o, \alpha u_1, \alpha u_2, w$ are the
vertices of a rhombus. The rotation $\psi$ rotates this rhombus onto
itself, exchanging the vectors $\alpha u_1$ and $\alpha u_2$.  The
points $o, \beta u_3, \beta u_4, -w$ form a rhombus through the same
axis, so that $\psi$ also exchanges the vectors $\beta u_3$ and $\beta
u_4$.  Since $\psi$ is a rotation, it preserves orthogonality.  It
follows that $P$ and $\psi P$ have the same normal vectors and the same
corresponding facet areas.  Proposition~\ref{special} then implies that $P$
and $\psi P$ are congruent by a translation.  In particular, the facets
$f_1$ and $f_2$ are congruent, as are $f_3$ and $f_4$. 
\end{proof}

The case of isosceles tetrahedra described in the introduction follows as an immediate
corollary to Theorem~\ref{pair}.
\begin{corollary} Suppose that $T$ is a tetrahedron in $\RR^3$.  
If the faces $f_i$ of $T$ satisfy the
condition
$$Area(f_1) \; = \; Area(f_2) \; = \; Area(f_3) \; = \; Area(f_4)$$
then $f_1 \cong f_2 \cong f_3 \cong f_4$.
\label{regular}
\end{corollary}
In other words, if a tetrahedron $T$ is equiareal, then $T$ is also
isosceles. For alternative proofs and variants 
of Corollary~\ref{regular}, see \cite{GemsII,Horv,Martini,equiareal}.

\noindent
{\bf Remark:}  Corollary~\ref{regular} has long been known to have an
analogue in which area is replaced by {\em perimeter}.  The proof is very
simple:  If all of the facets of $T$ have the same perimeter, 
the resulting system of linear equations (in the six edge lengths of 
$T$) implies that opposing edges must have the same length, so that
$T$ is isosceles.  A similar argument shows that
if the facets of $T$ can be partitioned into pairs
having the same perimeter then $T$ is reversible.

\section{\bf Factoring the volume}

Suppose that $T \subseteq \RR^3$ is a tetrahedron with vertices at 
$v_0, v_1, v_2, v_3 \in \RR^3$, 
where $v_0 = o$, the origin.
As before, let $A$ denote
the matrix whose columns are given by the vectors $v_i$, and suppose
that the $v_i$ are ordered so that $A$ has positive determinant.  
The volume of $T$ is then given by
$\det(A) = 6 V(T),$ so that
$$V(T)^2 = \frac{1}{36} \,\det(A^tA).$$

The entries of the matrix $A^tA$ are dot products of the form
$v_i \cdot v_j$.  From the identity,
\begin{equation}
2v_i \cdot v_j \;=\; |v_i|^2 + |v_j|^2 - |v_i-v_j|^2
\label{sq}
\end{equation}
it then follows that the value of $V(T)^2$ is a {\em polynomial} in
the {\em squares} of the edge lengths of $T$.  Said differently,
if $T$ has edge lengths $a_{ij}$ (the distance between vertices
$v_i$ and $v_j$),
then $V(T)^2$ is a polynomial in the variables $b_{ij} = a_{ij}^2$,
as well as the variables $a_{ij}$ themselves.
This polynomial is sometimes formulated in terms of linear algebraic
expressions such as Cayley-Menger determinants \cite[p.~125]{Sommer}.
While the Cayley-Menger heuristic outlined above applies in arbitrary
dimension, the $3$-dimensional case has been known at least as far back as 
Piero della Francesca \cite{Piero}.\footnote{Piero della Francesca (1415-1492), 
an Italian painter and geometer of the early Renaissance period.}

In certain instances, the polynomial $V(T)^2$ admits factorization into
linear or quadratic irreducible factors. For the 2-dimensional case, the
area $A(\Delta)$ of a triangle $\Delta$ having edge lengths $a,b,c$ is given by 
$$A(\Delta)^2 = \frac{1}{16}(a+b+c)(-a+b+c)(a-b+c)(a+b-c),$$
a factorization known as
{\em Heron's formula} \cite[p.~58]{Revisit}.  
Although the 3-dimensional case is more complicated \cite{DS},
there exist non-trivial factorizations of $V(T)^2$ when the tetrahedron $T$
satisfies the symmetry properties examined in the previous section.

For example, if $T$ is an isosceles tetrahedron, having edge lengths $a,b,c$
(each repeated twice in pairs of opposing edges), then
\begin{equation}
V(T)^2 = \frac{1}{72} \big( a^2 + b^2 - c^2 \big) 
\big( a^2 - b^2 + c^2 \big) \big( -a^2 + b^2 + c^2 \big).
\label{isovol}
\end{equation}
A synthetic proof of~(\ref{isovol}) can be found in \cite[p.~101]{Steinhaus}. 
Instead we will give an algebraic proof of the following more general
result, using a technique outlined in \cite{Klain-Heron}.

The edges of a reversible tetrahedron $T$ come in (at most) 4 lengths.  To see this,
label the edge lengths of $T$ so that the triangular facets 
$f_1 \cong f_2$ have edge lengths $a,b,c$, with common edge of length $c$. 
Since $f_3 \cong f_4$, they must have edge lengths $a,b,d$.
The six edges of $T$ then have lengths $a,a,b,b,c,d$, as in Figure~\ref{pn-fig}.

\begin{theorem}[Volume Formula] 
Suppose that $T$ is a reversible tetrahedron having 
edge lengths $a,a,b,b,c,d$.  Then
\begin{equation}
V(T)^2 
= \frac{1}{72} \bigg( c^2d^2 -  (a^2 - b^2)^2 \bigg) \bigg( a^2 + b^2 - \frac{c^2 + d^2}{2} \bigg).
\label{bivol}
\end{equation}
\label{pairvol}
\end{theorem}
The first polynomial factor in the formula~(\ref{bivol}) is a difference 
of two squares, so that~(\ref{bivol}) can be reformulated as
\begin{equation}
V(T)^2 = \frac{1}{72} \big( cd + a^2 - b^2 \big) \big( cd - a^2 + b^2 \big)
\bigg( a^2 + b^2 - \frac{c^2 + d^2}{2} \bigg).
\label{bivol2}
\end{equation}
In the special case where $c=d$, the tetrahedron $T$ is isosceles,
and the formula~({\ref{bivol2}) reduces to~(\ref{isovol}).

The proof of~(\ref{bivol}) will make use of two identities from plane geometry.
The well-known {\em parallelogram law} asserts that 
if edges of a parallelogram in $\RR^2$ are labelled as
in Figure~\ref{pn-fig}, then $2a^2+2b^2=c^2+d^2$.

The less well-known {\em trapezoid law}
asserts that, if the edges of a convex isosceles trapezoid are labelled as
in Figure~\ref{trap-fig}, then
$b^2-a^2 = cd$.

\begin{figure}[ht]
\begin{tikzpicture}[scale=1.2]
\draw[thick] (8,0)--(10,0)--(11,1.5)--(7,1.5)--(8,0);
\draw[thick] (8,0)--(11,1.5);
\draw[thick] (10,0)--(7,1.5);
\node [below] at (9,0.05) {\small $c$};
\node [above] at (9,1.45) {\small $d$};
\node [left] at (7.5,0.75) {\small $a$};
\node [right] at (10.5,0.75) {\small $a$};
\node [above] at (8.3,0.8) {\small $b$};
\node [above] at (9.7,0.8) {\small $b$};
\node [left] at (8,0) {\small $\oo$};
\node [right] at (10,0) {\small $\ww$};
\node [right] at (11,1.5) {\small $\vv$};
\node [left] at (7,1.5) {\small $\uu$};
\draw[fill] (8,0) circle [radius=0.05];
\draw[fill] (10,0) circle [radius=0.05];
\draw[fill] (11,1.5) circle [radius=0.05];
\draw[fill] (7,1.5) circle [radius=0.05];
\end{tikzpicture}
\caption{The trapezoid law: $b^2-a^2 = cd$.} 
\label{trap-fig}
\end{figure}
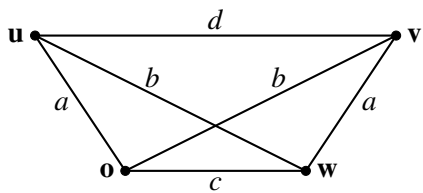
To see why, observe that
\begin{align*}
b^2-a^2 
&= |\uu-\ww|^2 - |\vv-\ww|^2\\ 
&= \uu \cdot \uu - 2\uu \cdot \ww + \ww \cdot \ww - \vv \cdot \vv + 2\vv \cdot \ww - \ww \cdot \ww\\ 
&= |\uu|^2  - |\vv|^2 + 2\ww \cdot (\vv-\uu)\\
&= a^2 - b^2 + 2cd, 
\end{align*}
where the last step follows from the parallelism of $\ww$ and $\vv-\uu$. 
The trapezoid law now follows.

\begin{proof}[Proof of The Volume Formula~\ref{pairvol}]
Let $f(a,b,c,d)$ denote the polynomial $V(T)^2$.  
The factors of $f$ can be determined by considering the cases in which
the volume of $T$ is zero, namely, when the tetrahedron $T$ is flat or 
otherwise degenerate.  If $T$ is reversible, this can occur in two ways.

In one case, $T$ may flatten to a parallelogram, having edges of length
$a,b,a,b$ and diagonals of length $c,d$.  In this instance, the
parallelogram law for the standard inner product implies that 
$2a^2 + 2b^2 = c^2 + d^2$.

In the second case, $T$ may flatten to a trapezoid, having non-parallel
edges of length $a,a$, parallel edges of length $c,d$, and diagonals of
length $b, b$.  In this instance, the 
trapezoid law implies that $(b^2 - a^2)^2 = c^2 d^2$,

These cases suggest both
$2a^2 + 2b^2 - c^2 - d^2$ and $c^2 d^2 - (b^2 - a^2)^2$
as possible factors of the polynomial $f$.

Denote $A = a^2$, $B = b^2$, $C = c^2$ and $D = d^2$.  
We observed following~(\ref{sq}) above that $f$ is 
a polynomial in the {\em squared} values $a^2,b^2,c^2,d^2$, so that
$f = f(A,B,C,D) \in \RR[A,B,C,D]$.  Since volume $V$ is homogeneous
of degree 3 with respect to length, the polynomial $f = V^2$
is homogeneous of degree 6 with respect to the variables
$a,b,c,d,$ and is therefore homogeneous of degree 3 with respect
to the variables $A,B,C,D$; that is, 
a homogeneous {\em cubic} polynomial in $\RR[A,B,C,D]$.

To verify that $2a^2 + 2b^2 - c^2 - d^2$ is indeed a factor of $f(a,b,c,d)$,
use division with remainder in $\RR[A,B,C,D]$ to obtain
$$f(A,B,C,D) = (2A+2B-C-D) g(A,B,C,D) + r(B,C,D),$$ 
for some $g \in \RR[A,B,C,D]$ and $r \in \RR[B,C,D]$.  
Here division with remainder in $\RR[A,B,C,D]$ is performed
here using lexicographical order on the variables 
$A,B,C,D$.  (See, for example, \cite[p.~54]{CLO}.)

Note that $A$ does not appear in the polynomial expression for $r$.  Suppose that
$C > D > 0$.  By the triangle inequality,
each $B$ such that 
$$\sqrt{C}-\sqrt{D} < 2\sqrt{B} < \sqrt{C}+\sqrt{D}$$ 
gives rise to a parallelogram as in Figure~\ref{para-fig},
yielding $A \geq 0$ so that $2A+2B-C-D=0$.  This
degenerate reversible
tetrahedron $T$ has volume zero, so that $f(A,B,C,D) = V^2 = 0$.  It
follows that $r(B,C,D) = 0$ on a non-empty open set.  Since $r$
is a polynomial, it follows that $r$ is identically zero, so that
$$f(A,B,C,D) = (2A+2B-C-D) g(A,B,C,D).$$
In other words, $2A+2B-C-D$ divides $f$ in $\RR[A,B,C,D]$.  

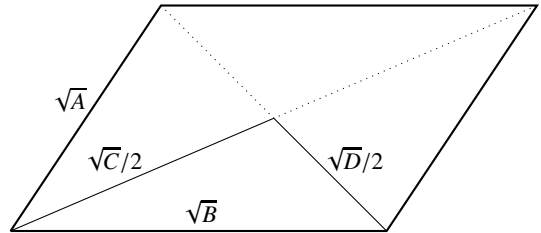
\begin{figure}[ht]
\begin{tikzpicture}[scale=1]
\draw[thick] (0,0)--(2,3)--(7,3)--(5,0)--(0,0);
\draw[thin] (0,0)--(3.5,1.5)--(5,0);
\draw[dotted] (2,3)--(3.5,1.5)--(7,3);
\node [left] at (1.15,1.75) {\scriptsize $\sqrt{A}$};
\node at (2.5,0.25) {\scriptsize $\sqrt{B}$};
\node [left] at (1.9,0.95) {\scriptsize $\sqrt{C}/2$};
\node [right] at (4,0.95) {\scriptsize $\sqrt{D}/2$};
\end{tikzpicture}
\caption{This parallelogram exists iff $\frac{1}{2}\sqrt{C}-\frac{1}{2}\sqrt{D} \leq \sqrt{B} \leq \frac{1}{2}\sqrt{C}+\frac{1}{2}\sqrt{D}$.} 
\label{para-fig}
\end{figure}

For the trapezoidal factors, view $f$ as polynomial in  $\RR[A,B,c,d]$, and write
$$f(A,B,c,d) = (cd-B+A) \tilde{g}(A,B,c,d) + \tilde{r}(B,c,d),$$
using division with remainder in $\RR[A,B,c,d]$ 
under lexicographical order on the variables 
$A,B,c,d$.  Once again the remainder $\tilde{r}$ is independent of the variable $A$,
while a trapezoidal degenerate (zero volume) tetrahedron can be constructed for an open set of
values $(B,c,d)$, so that $\tilde{r}$ is also identically zero. Therefore, $(cd-B+A)$ is also a factor $f$.

Finally, a symmetrical argument (reversing the roles of $A$ and $B$) yields a factor of $(cd-A+B)$.

Since $\RR[A,B,c,d]$ is a unique factorization domain \cite[p.~371]{Artin}\cite[p.~149]{CLO},
the irreducible factors 
$(cd-B+A)$, and $(cd-A+B)$ are prime, so that 
$$(cd-B+A)(cd-A+B) \;=\; c^2d^2-(B-A)^2 \;=\; CD-(B-A)^2$$
divides $f$.

Similarly, since $\RR[A,B,C,D]$ is a unique factorization domain,
the two irreducible factors 
$CD-(B-A)^2$, and $2A+2B-C-D$ are prime in $\RR[A,B,C,D]$, 
so that 
\begin{align}
V^2 \;=\; f \;=\; (2A+2B-C-D)(CD - (B-A)^2)k.
\label{um}
\end{align}
Because $f$ is a homogeneous cubic polynomial in $\RR[A,B,C,D]$, the factor $k$ must be a constant,
independent of the parameters $A,B,C,D$.

To compute the constant $k$, recall that the volume of the regular
(equilateral) tetrahedron of unit edge length 
$A = B = C = D = 1$ is $\sqrt{2}/12$.  It
follows that 
$$\frac{1}{72} = \left(\frac{\sqrt{2}}{12}\right)^2 = 
V^2 = f(1,1,1,1) = 2k.$$
Hence, $k = 1/144$, and~(\ref{um}) becomes~(\ref{bivol}).
\end{proof}

I.~Izmestiev has pointed out that
applying the {\em Regge symmetry} \cite{ako1} to a reversible tetrahedron 
gives a new reversible tetrahedron having the same volume, 
and for which the factors of the Cayley-Menger polynomial~(\ref{bivol2}) are permuted \cite{II}.

\section{\bf Generalizations}
\label{general}

A convex polytope $P$ in $\RR^n$ will be called 
{\em reversible}
if there is an affine plane $\xi$ of co-dimension $2$ such that
$P$ is symmetric under the $180^\circ$ rotation of the $2$-plane $\xi^\perp$ that fixes $\xi$.  

If a tetrahedron $T$ is
$\RR^3$ is symmetric under a $180^\circ$ rotation around a line $\ell$, 
then this rotation must map facets to facets and facet normals to facet normals.
In view of Proposition~\ref{special}, 
the only way this can occur is when $\ell$ passes through the midpoints of two non-adjacent edges of $T$, 
so that $T$ must have pairs of congruent facets, as in the examples addressed earlier.  
It follows that this more general definition of a reversible polytope is
consistent with the definition given earlier for tetrahedra in $\RR^3$.  
However, naive analogues of the theorems of this paper do not follow, 
because this level of symmetry admits many more variations in structure for dimensions $n \geq 4$.  
Indeed, there exist $4$-dimensional simplices in which all $5$ facets have the same volume in spite of not being mutually congruent.  
For an extensive treatment of this subject, see \cite{equiareal}.

In addition to admitting the Heron-type formula~(\ref{isovol}) for volume, isosceles tetrahedra satisfy many other characteristic properties 
(see, for example, \cite[p. 90-97]{GemsII}\cite{Leech}).  
It would be interesting to consider what parallels these other properties may have in the more general context of reversible tetrahedra.


\begin{thebibliography}{1}

\bibitem{ako1} A.~Akopyan and I.~Izmestiev,
The Regge symmetry, confocal conics, and the Schl\"{a}fli formula,
{\em Bull. London Math. Soc.}, {\bf 51} (2019), 765--775.


\bibitem{Artin}
M.~Artin, \emph{Algebra, 2nd ed.}, Prentice-Hall, Upper Saddle River, NJ, 2010.

\bibitem{Bonn}
T.~Bonnesen and W.~Fenchel, 
\emph{Theory of {C}onvex {B}odies}, BCS Associates, Moscow, Idaho, 1987.

\bibitem{Coxeter}
H.~Coxeter, \emph{Regular {P}olytopes}, Dover, New York, 1973.

\bibitem{Revisit}
H.~Coxeter and S.~Greitzer, \emph{Geometry {R}evisited}, MAA, Washington, D.C., 1967.

\bibitem{CLO} D.~Cox, J.~Little, and D.~O'Shea,  \emph{Ideals, {V}arieties,
and {A}lgorithms, 2nd ed.}, Springer, New York, 1996.

\bibitem{DS} C.~D'Andrea and M.~Sombra,
The Cayley-Menger determinant is irreducible for $n \geq 3$,
{\em Sib.~Math.~J.}, {\bf 46} (2005), 71--76.

\bibitem{FR1} P.~Fowler and A.~Rasset, Is There a ``Most Chiral Tetrahedron''?,
{\em Chem. Eur. J.}, {\bf 10} (2004), 6575--6580.

\bibitem{FR2} P.~Fowler and A.~Rasset, A classification scheme for chiral tetrahedra,
{\em C. R. Chimie}, {\bf 9} (2006), 1203--1208.

\bibitem{GemsII}
R.~Honsberger, \emph{Mathematical {G}ems {I}{I}}, MAA, Washington, D.C., 1976.

\bibitem{Horv}
J. ~Horv\'ath, A property of tetrahedra with equal faces in spaces of constant curvature (Hungarian,
German summary), {\em Mat. Lapok}, {\bf 20} (1969), 257--263.

\bibitem{II} I.~Izmestiev. {\em Private communication}, (2022). 

\bibitem{Klain-Heron} D.~Klain, 
An intuitive derivation of Heron's formula,
{\em Amer. Math. Monthly}, {\bf 111} (2004), no. 8, 709--712. 

\bibitem{Leech} J.~Leech, 
Some properties of the isosceles tetrahedron,
{\em Math. Gaz.}, {\bf 34} (1950), no. 310, 269--271. 

\bibitem{Martini} H.~Martini,
Regular simplices in spaces of constant curvature,
{\em Amer. Math. Monthly}, {\bf 100} (1993), no. 2, 169--171. 

\bibitem{equiareal} P.~McMullen, Simplices with equiareal faces,
{\em Discrete Comput. Geom.}, {\bf 24} (2000), 397--411.

\bibitem{Piero} M.~Peterson, The geometry of Piero della Francesca, 
{\em Math. Intell.}, {\bf 19} (1997), no. 3, 33--40. 

\bibitem{red}
R.~Schneider, \emph{Convex {B}odies: {T}he {B}runn-{M}inkowski {T}heory, 2nd ed.},
  Cambridge University Press, New York, 2014.

\bibitem{Sommer}  D.~Sommerville,
\emph{An {I}ntroduction to the {G}eometry of n {D}imensions},
Dover, New York, 1958. 

\bibitem{Steinhaus}
H. Steinhaus,
\emph{One {H}undred {P}roblems in {E}lementary {M}athematics},
Dover, New York, 1979.

\bibitem{WD1} K.~Wirth and A.~Dreiding,
Edge lengths determining tetrahedrons, {\em Elem.~Math.}, {\bf 64} (2009), 160--170.

\bibitem{WD2} K.~Wirth and A.~Dreiding,
Tetrahedron classes based on edge lengths, {\em Elem.~Math.}, {\bf 68} (2013), 56--64.

\bibitem{WD3} K.~Wirth and A.~Dreiding,
Relations between edge lengths, dihedral, and solid angles
in tetrahedra, {\em J.~Math.~Chem.}, {\bf 52} (2014), 1624 --1638.

\end{thebibliography}
\end{document}